\documentclass[a4paper]{article} 
\usepackage{graphicx} 

\usepackage[russian]{babel}

\long\def\comment#1{}

\makeatletter
\def\@normalsize{\@setsize\normalsize{10pt}\xpt\@xpt
\abovedisplayskip 10pt plus2pt minus5pt\belowdisplayskip
\abovedisplayskip \abovedisplayshortskip \z@
plus3pt\belowdisplayshortskip 6pt plus3pt
minus3pt\let\@listi\@listI}

\def\subsize{\@setsize\subsize{12pt}\xipt\@xipt}
\def\section{\@startsection {section}{1}{\z@}{1.0ex plus
1ex minus .2ex}{.2ex plus .2ex}{\large\bf}}
\def\subsection{\@startsection
   {subsection}{2}{\z@}{.2ex plus 1ex} {.2ex plus .2ex}{\subsize\bf}}
\makeatother

\begin{document}

\date{}

\title{\huge \bf {Индуцированные структуры на трансверсалях к подгруппе группы
и гипергруппы над группой}}


\author{\bf{Самвел Г. Далалян}\\
Ереванский государственный университет\\\\
{\bf Аннотация.} На трансверсалях к подгруппе группы  с помощью бинарной операции\\
группы определяются структурные отображения и на их основе вводится понятие\\
гипергруппы над группой, обобщающее понятие фактор-группы заданной группы по\\
 нормальной подгруппе. Категория гипергрупп над группой содержит в качестве\\
подкатегорий категорию групп, категорию линейных пространств, категорию полей (тел).\\
{\bf Ключевые слова}: трансверсаль, фактор-группа, гипергруппа над группой, линейное\\
пространство, поле.
}
  \maketitle
	
	УДК 512548 +512538


\hspace{10mm}
0. Одной из фундаментальных конструкций теории групп является
конструкция фактор-группы данной группы по ее нормальной подгруппе.
По любой данной группе $G$ и произвольной ее нормальной подгруппе $H$
строится фактор-множество $G \slash H$, затем оно снабжается индуцированной структурой
группы.
Фактор-множество $G \slash H$ интерпретируется как разбиение множества $G$ на множество
смежных классов группы $G$ по нормальной подгруппе $H$, произвольным образом фиксируется
по одному элементу в каждом смежном классе и определяется бинарная операция на множестве
$G \slash H$ с помощью бинарной операции  группы $G$.

\hspace{10mm}
При этом существенную роль играют следующие два обстоятельства.
Во-первых, по произвольной, в том числе нормальной, подгруппе $H$ группы $G$ и любому элементу
$a \in G$ можно построить два смежных класса: правый $Ha = \{ x \cdot a, x \in H \}$ и левый
$aH = \{ a \cdot x, x \in H \}$; в случае нормальной подгруппы $H$ эти два смежных класса совпадают.
Во-вторых,   произведение смежных классов по вышеуказанному принципу в случае нормальной
подгруппы $H$  определяется корректно, то есть не зависит от выбора   отмеченных элементов.

\hspace{10mm}
Если отказаться от условия нормальности подгруппы $H$, эти два свойства уже не выполняются.
Оказывается, тем не менее, на подходяще выбранных множествах $M$ можно задать бинарную операцию,
индуцированную бинарной операцией группы $G$,
однако, в общем случае, эта операция не будет групповой.
Данная заметка посвящена конструкции этой бинарной операции и трех сопутствующих структурных отображений,
на основе которых вводится понятие гипергруппы над группой.
Доказывается, что категория гипергрупп над
группой содержит в качестве подкатегорий категорию групп, категорию линейных пространств, категорию
полей (тел).

\vspace{5mm}


\hspace{10mm}
{\bf 1}.
Подходящий для наших целей выбор множества $M$ основан на понятии дополнительного множества к подгруппе группы.
Пусть $G$ - группа, $H$ - ее произвольная подгруппа.
Подмножество $M$ группы $G$ называется
{\it правым дополнительным множеством к подгруппе $H$}, если
для любого элемента $x \in G$ существуют единственные элементы $\alpha \in H$ и $a \in M$,
такие что $x = \alpha \cdot a$.
Аналогично (двойственно) определяется понятие левого дополнительного множества к подгруппе $H$.
Вообще,  основные понятия и конструкции настоящей заметки имеют два варианта: правый и левый.
В данной заметке  мы, как правило, рассматриваем только соответствующие правые варианты и
слово правый часто опускаем.

\hspace{10mm}
Правой {\it трансверсалью к подгруппе $H$ группы $G$} называется подмножество $M$ группы $G$,
пересекающее каждый правый смежный класс относительно подгруппы $H$ в точности по одному элементу.

\hspace{10mm}
{\it Сечением} канонического (правого) сюръективного фактор-отображения
$$
\psi: G \rightarrow H \backslash G, \, a \mapsto Ha
$$
называется отображение $\sigma: H \backslash G \rightarrow G$, удовлетворяющее условию
$\sigma \circ \psi = id$, где $id$ означает тождественное отображение.

\hspace{10mm}
Следующие три утверждения эквивалентны:

$M$ - правое дополнительное множество к подгруппе  $H$ группы $G$;

$M$ - правая трансверсаль  к подгруппе  $H$ группы $G$;

$M$ - образ сeчения $\sigma$ правого фактор-отображения $\psi$ (см. [1]).

Поэтому термины  правое дополнительное множество, правая трансверсаль
и образа сечения правого фактор-отображения  взаимозаменяемы.


\vspace{5mm}

\hspace{10mm}
{\bf 2.}
Пусть $G$ - произвольная группа, $H$ - ее подгруппа, $M$ - дополнительное множество  к подгруппе $H$.
Мы обозначаем элементы группы $H$ строчными греческими буквами, элементы множества $M$ -
строчными латинскими буквами.

\hspace{10mm}
Определим отображения
\begin{itemize}
\item[\bf{(B1)}] $\Phi: M \times H \rightarrow M, \quad \Phi (a, \alpha) : = a^\alpha$,

\item[\bf{(B2)}] $\Psi: M \times H \rightarrow H, \quad \Psi (a, \alpha) := {^a}\alpha$,

\item[\bf{(B3)}] $\Xi: M \times M \rightarrow M, \quad \Xi (a, b) := [a, b]$,

\item[\bf{(B4)}] $\Lambda: M \times M \rightarrow H, \quad \Lambda (a, b) := (a, b)$
\end{itemize}
соотношениями
$$
(F1) \hspace{40mm}
a \cdot \alpha = {^a}\alpha \cdot a^\alpha, \quad   a \cdot b = (a, b) \cdot [a, b].
\hspace{50mm}
$$


\vspace{5mm}

\hspace{10mm}
{\bf 3.}
{\bf Теорема.}
{\it Отображения $\Phi, \Psi, \Xi, \Lambda$ обладают следующими свойствами.}

(P1) {\it Отображение $\Xi$ определяет на $M$ структуру
правой квазигруппы с левым нейтральным элементом $o$.}

(P2) {\it Отображение $\Phi$ является правым действием группы $H$ на множество $M$.}

(P3) {\it Отображение $\Psi$ переводит подмножество $\{ o \} \times H$ на $H$.}

(P4) {\it Система отображений $\Omega = (\Phi, \Psi, \Xi, \Lambda)$ подчиняется следующим
соотношениям}:
\begin{eqnarray*}
(A1) &{^a}(\alpha \cdot \beta) = {^a}\alpha \cdot {^{a^\alpha}}\beta, \\
(A2) &[a,b]^{\alpha} = [a^{^{b}\alpha}, b^{\alpha}],\\
(A3) &(a,b) \cdot ^{[a,b]}\alpha = {^a}(^{b}\alpha)\cdot (a^{^{b}\alpha}, b^{\alpha}),\\
(A4) &[[a, b], c] = [a^{(b,c)}, [b, c]],\\
(A5) &(a, b) \cdot ([a, b], c) = {^a}(b, c) \cdot (a^{(b, c)}, [b, c]),\\
\end{eqnarray*}

\hspace{10mm}
{\it Доказательство.}
Используем  ассоциативность бинарной
операции группы $G$ в виде
$$
(a \cdot \alpha) \cdot \beta = a \cdot (\alpha \cdot \beta), \quad
a \in M, \, \alpha, \beta \in H.
$$
Применяя $(F1)$ к левой части, получаем
$$
(a \cdot \alpha) \cdot \beta = ({^a}\alpha \cdot a^{\alpha}) \cdot \beta =
{^a}\alpha \cdot (a^{\alpha} \cdot \beta) =
{^a}\alpha \cdot ({^{a^{\alpha}}}\beta \cdot (a^{\alpha})^{\beta}) =
({^a}\alpha \cdot {^{a^{\alpha}}}\beta) \cdot (a^{\alpha})^{\beta}.
$$
Применив $(F1)$ к правой части, получим
$$
a \cdot (\alpha \cdot \beta) = {^a}(\alpha \cdot \beta) \cdot a^{\alpha \cdot \beta}.
$$
Следовательно,
$$
({^a}\alpha \cdot {^{a^{\alpha}}}\beta) \cdot (a^{\alpha})^{\beta} = {^a}(\alpha \cdot \beta) \cdot a^{\alpha \cdot \beta}, \quad
{^a}\alpha \cdot {^{a^{\alpha}}}\beta, \, {^a}(\alpha \cdot \beta) \in H, \quad (a^{\alpha})^{\beta}, \, a^{\alpha \cdot \beta} \in M.
$$
Отсюда, используя единственность разложения в произведение элементов подгруппы $H$ и трансверсали $M$, получаем
соотношения $(A1)$ и
\begin{center}
$(A0) \hspace{15mm} (a^{\alpha})^{\beta} = a^{\alpha \cdot \beta}. \hspace{20mm}$
\end{center}
Соотношение  (А0) является первым необходимым соотношением определения действия группы  на множество (свойство (P2)).

\hspace{10mm}
Аналогично, используя ассоциативность бинарной операции группы $G$ в частных случаях
$$
(a \cdot b) \cdot \alpha = a \cdot (b \cdot \alpha), \quad a, b \in M, \,\, \alpha \in H
$$
и
$$
(a \cdot b) \cdot c = a \cdot (b \cdot c), \quad a, b, c \in M,
$$
получим соотношения $(A2), (A3)$ и $(A4), (A5)$, соответственно.

\hspace{10mm}
Применяя те же идеи доказательства к равенству
\begin{center}
$a \cdot \varepsilon = \varepsilon \cdot a, \quad a \in M$, $\varepsilon$ -
нейтральный элемент подгруппы $H$,
\end{center}
можно показать, что
$a^{\varepsilon} = a, \, \, {^a}\varepsilon = \varepsilon.$
Первое из этих соотношений комплектует свойство (P2).

\hspace{10mm}
Осталось доказать свойства (P1) и (P3).
Пусть
$$
e = \theta \cdot o, \quad \theta \in H, o \in M
$$
каноническое разложение нейтрального элемента $e = \varepsilon$ группы $G$
в произведение элементов подгруппы $H$ и трансверсали $M$.
Тогда для любого элемента $a \in M$ имеем
$$
\varepsilon \cdot a = \theta \cdot o \cdot a = (\theta \cdot (o, \, a)) \cdot[o, \, a].
$$
Следовательно, в частности, $[o, a] = a$.
Значит бинарная операция $\Xi$ имеет левый нейтральный элемент $o$.

\hspace{10mm}
Так как
$$
(\alpha \cdot \theta) \cdot o = \alpha = \theta \cdot (o \cdot \alpha) =
(\theta \cdot {^o}\alpha) \cdot o^\alpha,
$$
получаем соотношения ${^o}\alpha = \theta^{-1} \cdot \alpha \cdot \theta$
и $o^\alpha = o$.
Свойство (P3) является следствием первого из этих соотношений.

\hspace{10mm}
Наконец, докажем, что множество $M$ вместе с бинарной операцией $\Xi$
является правой квазигруппой, то есть любое уравнение $[x, a] = b, \,\, a, b \in M$
имеет единственное решение в $M$.

\hspace{10mm}
Для любого элемента $a \in M$ элементы $a^{(-1)} \in H$ и $a^{[-1]} \in M$ определим условием
$$
a^{-1} = a^{(-1)} \cdot a^{[-1]}.
$$
Тогда
$$
a^{(-1)} \cdot (a^{[-1]}, \, a) \cdot [a^{[-1]}, \, a] = a^{(-1)} \cdot a^{[-1]} \cdot a
= a^{-1} \cdot a = e = \theta \cdot o,
$$
и, следовательно,
$$
(F2) \hspace{20mm} [a^{[-1]}, \, a] = o, \quad a^{(-1)} \cdot (a^{[-1]}, \, a) = \theta. \hspace{30mm}
$$

\hspace{10mm}
{\bf Лемма.}
{\it Для  элементов $x, a, b \in M$  соотношение $[x, a] = b$ выполняется
тогда и только тогда, когда
$$
x = [b^{a^{(-1)}}, a^{[-1]}], \quad
(x, a) \cdot {^b}(a^{(-1)}) \cdot (b^{a^{(-1)}}, a^{[-1]})  = \varepsilon.
$$}

\hspace{10mm}
{\it Доказательство леммы.}
На основании (F1) и единственности разложения в произведение элементов подгруппы $H$
и трансверсали к ней $M$ имеем цепочку эквивалентных равенств:
$$
[x, a] = b \quad \Leftrightarrow \quad
(x, a) \cdot [x, a] = (x, a) \cdot b \quad \Leftrightarrow \quad
x \cdot a = (x, a) \cdot b \quad \Leftrightarrow \quad
$$
$$
x = (x, a) \cdot b \cdot a^{-1} \quad \Leftrightarrow \quad
x = (x, a) \cdot b \cdot a^{(-1)} \cdot a^{[-1]} \quad \Leftrightarrow \quad
$$
$$
x = (x, a) \cdot {^b}(a^{(-1)}) \cdot b^{a^{(-1)}} \cdot a^{[-1]} \quad \Leftrightarrow \quad
$$
$$
\varepsilon \cdot x = (x, a) \cdot {^b}(a^{(-1)}) \cdot (b^{a^{(-1)}}, a^{[-1]}) \cdot [b^{a^{(-1)}}, a^{[-1]}] \quad \Leftrightarrow \quad
$$
$$
x = [b^{a^{(-1)}}, a^{[-1]}], \quad  (x, a) \cdot {^b}(a^{(-1)}) \cdot (b^{a^{(-1)}}, a^{[-1]}) = \varepsilon.
$$

\hspace{10mm}
{\bf Следствие.}
{\it Уравнение $[x, a] = b$, $a, b \in M$ может иметь не более одного решения в $M$.
Таким решением может быть только $x = [b^{a^{(-1)}}, a^{[-1]}]$.
Это значение неизвестного будет решением рассматриваемого уравнения тогда и
только тогда, когда имеет место равенство
$$
(F3) \hspace{20mm}
([b^{a^{(-1)}}, a^{[-1]}], a) \cdot {^b}(a^{(-1)}) \cdot (b^{a^{(-1)}}, a^{[-1]})
 = \varepsilon.
 \hspace{35mm}
$$}

Подставляя значение  $x = [b^{a^{(-1)}}, a^{[-1]}]$ в уравнение $[x, a] = b$, убеждаемся,
что оно является его решением:
$$
[[b^{a^{(-1)}}, a^{[-1]}], a] = [[b^{a^{(-1)} \cdot (a^{[-1]}, a)}, [a^{[-1]}, a]] =
[b^\theta, o] = (b^\theta)^{\theta^{-1}} = b.
$$
Здесь используются $(A4), (F2)$ и соотношение $[a, o] = a^{\theta^{-1}}$,
которое (вместе с соотношением $(a, o) = {^a}(\theta^{-1})$)
можно извлечь из равенств
$$
(a, o) \cdot [a, o] = a \cdot o = a \cdot \theta^{-1} = {^a}(\theta^{-1}) \cdot a^{\theta^{-1}}.
$$
Таким образом, теорема полностью доказана.
Заметим, что, как следствие,  справедливо соотношение $(F3)$.


\vspace{5mm}

\hspace{10mm}
{\bf 4.}
{\bf Определение}  (см. [2], [4]).
Пара $(M, H)$, состоящая из множества $M$ и группы $H$,
вместе с системой структурных отображений $\Omega = (\Phi, \Psi, \Xi, \Lambda)$
называется {\it правой гипергруппой над группой $H$}
(и обозначается $M_H$),
если выполняются свойства (P1) - (P4).

\hspace{10mm}
Теперь нашу основную теорему 3 можно переформулировать так.

\hspace{10mm}
{\it С каждой тройкой $(G, H, M)$, где $G$ - группа, $H$ - ее подгруппа,
$M$ - правая трансверсаль к подгруппе $H$, канонически ассоциируется правая
гипергруппа над группой $M_H$.
}

\hspace{10mm}
Конструкция, с помощью которой мы получили гипергруппу над группой $M_H$,
исходя из тройки $(G, H. M)$,
называется {\it стандартной конструкцией} гипергруппы над группой.
Можно доказать (см, [6]), что стандартная конструкция {\it универсальна}:
с точностью до изоморфизма, любая гипергруппа над группой
получается с помощью стандартной конструкции
(понятие изоморфизма гипергрупп над группой определяется в  п. 7
настоящей заметки).

\hspace{10mm} Термин гипергруппа уже используется в математике.
Так называются объекты со структурой, похожей на групповую,
но с многозначной бинарной операцией
(см., например, [7]; о некоторых применениях гипергрупп см. в [8]).
Выбор термина для введенного нами объекта обосновывается тем, что
этот объект является обобщением группы и линейного пространства
над полем (телом) (см. п. 9 и 10 этой заметки).
Путаницы с понятием гипергруппы не может возникнуть,
потому что (полное) название введенного объекта - {\it гипергруппа над группой}.


\vspace{5mm}

\hspace{10mm}
{\bf 5.}
Сделаем несколько комментариев к соотношениям $(A1) - (A5)$ определения гипергруппы над группой.

\hspace{10mm}
(a1) Условие $(A1)$ означает, что для любого элемента $a \in M$ отображение
$$
\Psi_a: H \rightarrow H, \quad \Psi_a (\alpha) = {^a}\alpha
$$
является "`обобщенным гомоморфизмом"'.
Оно превращается в обычный гомоморфизм, если
$\Phi_{\alpha} (a) = a^\alpha = a$ при любом $\alpha \in H$,
то есть когда  $a$ - неподвижный элемент относительно действия $\Phi$ группы $H$ на множество $M$.

\hspace{10mm}
(a2) Аналогично условие $(A2)$ означает, что для любого элемента $\alpha \in H$ отображение
$$
\Phi_{\alpha}: M \rightarrow M, \quad \Phi_{\alpha} (a) = a^{\alpha}
$$
является "`обобщенным гомоморфизмом"' для бинарной  операции $\Xi (a, b) = [a, b]$ квазигруппы $(M, \Xi)$.
Оно превращается в обычный гомоморфизм, если ${^{b}}\alpha = \alpha$ для любого $b \in M$.

\hspace{10mm}
(a3) Формула $(A3)$ является аналогом формулы
$(a^{\alpha})^{\beta} = a^{\alpha \cdot \beta}$.
Кроме того имеем аналог
$$
(F4) \hspace{30mm} {^o}\alpha = \theta^{-1} \cdot \alpha \cdot \theta  \hspace{35mm}
$$
формулы $a^{\varepsilon} = a$.
Значит отображение $\Psi$ можно рассматривать как некое обобщенное кодействие
правой квазигруппы $(M, \Xi)$ с левым нейтральным элементом $o$
на группу $H$ обобщенными автоморфизмами.

\hspace{10mm}
Отметим, что при доказательстве теоремы 3 мы использовали соотношение  $(F4)$
 для получения свойства $(P3)$.
На самом деле, здесь справедлива и обратная импликация, так что соотношение $(F4)$
 и свойство $(P3)$ эквивалентны.
Условие $(P3)$ в определении гипергруппы над группой просто более удобная форма
задания соотношения $(F4)$.

\hspace{10mm}
(a4) Соотношение $(A4)$, очевидно,  формула обобщенной ассоциативности бинарной операции $\Xi$.

\hspace{10mm}
(a5) Вскользь отметим, что соотношение $(A5)$ можно рассматривать как обобщение коцепного правила
теории когомологий групп (см. [9], [1], [6]).


\vspace{5mm}

\hspace{10mm}
{\bf 6.}
{\bf Предложение.}
{\it Пусть $H$ - нормальная подгруппа группы $G$,
$M$ - произвольная трансверсаль к подгруппе $H$,
$\Omega = (\Phi, \Psi, \Xi, \Lambda)$ - система
структурных отображений гипергруппы над группой $M_H$.
Тогда

1)  действие $\Phi$ группы $H$ на множество $M$ - тривиальное,
то есть $a^{\alpha} = a$ для любых $a \in M, \alpha \in H$;

2) правая квазигруппа $(M, \Xi)$ является группой,
изоморфной фактор-группе $H \backslash G$.}

\hspace{10mm}
{\it Доказательство.}
1). По определению нормальной подгруппы для любых элементов
$\alpha \in H, a \in M$ элемент $a \cdot \alpha \cdot a^{-1}$ должен принадлежать $H$.
Значит $a^{\alpha} = a$.

\hspace{10mm}
2). Если действие $\Phi$ - тривиальное, то согласно $(A4)$
правая квазигруппа $(M, \Xi)$ с левым нейтральным элементом ассоциативна.
Но любая ассоциативная правая квазигруппа с левым нейтральным элементом является группой.

\hspace{10mm}
Фактор-отображение $\psi: G \rightarrow H \backslash G, \, x \mapsto H \cdot x$
устанавливает биекцию между множествами $M$ и $H \backslash G$ и переводит бинарную операцию $\Xi$
в бинарную операцию фактор-группы $H \backslash G$.
Значит группа $(M, \Xi)$ изоморфна фактор-группе $H \backslash G$,
причем это верно для любой трансверсали $M$ к подгруппе $H$.
Предложение 5 доказано.

\hspace{10mm}
6.1. {\bf Следствие.}
{\it Пусть $H$ - нормальная подгруппа группы $G$.
Пусть $M'$ и $M''$, соответственно, правая и левая трансверсали к подгруппе $H$, а
$\Xi'$ и $\Xi''$ - соответствующие бинарные операции.
Тогда группы $(M', \Xi')$ и $(M'', \Xi'')$ изоморфны.}

\hspace{10mm}
{\it Доказательство.} Группа $(M', \Xi')$ изоморфна фактор-группе $H \backslash G$,
группа $(M'', \Xi'')$ изоморфна фактор-группе $G \slash H$, а фактор-группы
$H \backslash G$  и  $G \slash H$ совпадают.

\hspace{10mm}
Итак, понятие гипергруппы над группой является обобщением понятия фактор-группы по нормальной подгруппе,
обобщением в смысле предложения 6.

\hspace{10mm}
6.2. {\bf Следствие.} {\it Для всякой тройки $(G, H, M)$, где $H$ - подгруппа группы $G$  индекса 2, квазигруппа
$(M, \Xi)$ соответствующей  этой тройке гипергруппы над группой $M_H$ будет группой.}

\hspace{10mm}
{\it Доказательство.} Всякая подгруппа индекса 2 - нормальная.

{\hspace{10mm} }
Таким образом, согласно свойству универсальности стандартной конструкции, гипергруппы над группой с
некоммутативной (эквивалентно, не групповой) квазигруппой $(M, \Xi)$ могут возникнуть, только если
число элементов $|M|$ множества $M$ больше двух.
Все гипергруппы над группой $M_H$ с $|M| = 3$, с  точностью до изоморфизма, описаны в кандидатской
диссертации П. Зольфагари.
Некоторые примеры гипергрупп над группой $M_H$ с $|M| = 3$ и некоммутативной квазигруппой $(M, \Xi)$
приведены в [5].


\vspace{5mm}

\hspace{10mm}
{\bf 7.}
Класс гипергрупп над группой можно превратить в категорию,
определив понятие морфизма гипергрупп над группой (см. [3], [6]).

\hspace{10mm}
Пусть $M_H$ и $M'_{H'}$ две гипергруппы над группой с системами
структурных отображений $\Omega = (\Phi, \Psi, \Xi, \Lambda)$
и $\Omega' = (\Phi', \Psi', \Xi', \Lambda')$, соответственно.
{\it Морфизм из $M_H$ в $M'_{H'}$} - это пара $f = (f_0, f_1)$,
состоящая из гомоморфизма групп $f_0: H \rightarrow H'$ и
отображения множеств $f_1: M \rightarrow M'$, которые подчиняются
соотношениям
\begin{center}
$\Phi \cdot f_{1}$ = $(f_{1}\times f_{0})\cdot \Phi^{'}$, \quad
$\Psi \cdot f_{0}$ = $(f_{1}\times f_{0}) \cdot \Psi^{'}$, \quad
$\Xi \cdot f_1 = (f_1 \times f_1) \cdot \Xi', \quad$
$ \Lambda \cdot f_{0}$ = $(f_{1}\times f_{1}) \cdot \Lambda^{'}$.
\end{center}
Классы всех гипергрупп над группой и их морфизмов
вместе с покомпонентной операцией композиции морфизмов
$f \circ g = (f_0 \circ g_0, f_1 \circ g_1)$, где
$f = (f_0, f_1), g = (g_0, g_1)$,
и системой тождественных морфизмов $1_{M_H} = (1_H, 1_M)$
образуют категорию, которую мы обозначим $\mathcal{H}g$.
Поэтому в отношении гипергрупп над группой
мы можем использовать всю терминологию и все результаты
общей теории категорий.

\hspace{10mm}
Обозначим категорию гипергрупп над группой через $\mathcal{H}g$.


\vspace{5mm}

\hspace{10mm}
{\bf 8.}
{\bf Предложение.}
{\it Категория гипергрупп над группой $\mathcal{H}g$ содержит подкатегорию,
изоморфную категории групп $\mathcal{G}$.
Это - полная подкатегория с классом объектов, состоящим из всех гипергрупп
$M_E$ над тривиальной группой $E$.}

\hspace{10mm}
{\it Доказательство.}
Определим функтор $S(G) : \mathcal{G} \rightarrow \mathcal{H}g$ следующим образом.
Каждой группе $M$ сопоставим пару $(M, E)$ с тривиальной (под)группой $E$
и соответствующую гипергруппу $M_E$.
Очевидно, что при этом структурные отображения $\Phi, \Psi, \Lambda$ определяются
однозначно (и тривиально), а в качестве $\Xi$ возьмем бинарную операцию группы $M$.
Тогда выполняются условия $(P1) - (P4)$, так что получаем гипергруппу над группой.
Гомоморфизму групп $f_1: M \rightarrow M'$ сопоставим морфизм гипергрупп над группой
$$
f = (1_E, f_1): M_E \rightarrow M'_E.
$$
Этот функтор инъективен, поэтому устанавливает изоморфизм между
категорией $\mathcal G$ и ее образом ${\mathcal H}g_E$ относительно функтора $S(G)$,
являющимся полной подкатегорией категории $\mathcal{H}g$,
классом объектов которой служат все гипергруппы над тривиальной
группой $E$. Обратным функтором к функтору $S(G)$ будет стирающий функтор из  (под)категории
${\mathcal H}g_E$, забывающий структурные отображения $\Phi, \Psi, \Lambda$ и тривиальную группу $E$.


\vspace{5mm}

 \hspace{10mm}
{\bf 9.}
{\bf Предложение.}
{\it Пусть $k$ - произвольное поле.
Категория гипергрупп над группой $\mathcal{H}g$ содержит подкатегорию,
изоморфную категории линейных пространств $\mathcal{L}_k$ над фиксированным
полем $k$.}

\hspace{10mm}
{\it Доказательство.}
Пусть $H$ - мультипликативная группа поля $k$.
Каждому линейному пространству $L$ над полем $k$ сопоставим гипергруппу
$M_H$ над группой $H$ с системой структурных отображений $\Omega = (\Phi, \Psi, \Xi, \Lambda)$,
где $M$ совпадает с базовым множеством  пространства $L$,
$\Phi (a, \alpha) = a \alpha$ - произведение вектора $a \in M = L$ на скаляр $\alpha \in  H \subset k$,
$\Xi$ - бинарная операция аддитивнoй группы  пространства $L$,
а $\Psi$ и $\Lambda$  тривиальны:
$$
   \Psi (a, \alpha) = \alpha, \quad \Lambda (a, b) = \varepsilon,
$$
где $\varepsilon = 1$ - единичный элемент поля $k$.
Из аксиом линейного пространства следует выполнение условий $(P1) - (P4)$
(часть из которых превращается в тавтологию).

\hspace{10mm}
Линейному отображению $f: L \rightarrow L'$ линейных пространств над полем $k$
сопоставим морфизм гипергрупп над группой $(1_H, f)$, где $1_H$ - тождественный гомоморфизм
группы $H$.
Построенная пара сопоставлений
$$
L_k \mapsto M_H, \quad f \mapsto (1_H, f)
$$
определяет инъективный функтор $S(L)$ из категории ${\mathcal L}_k$ в категорию ${\mathcal H}g$.
Его образ $im S(L)$ будет подкатегорией категории ${\mathcal H}g$, изоморфной категории ${\mathcal L}_k$.

\hspace{10mm}
Подкатегория $im S(L)$ имеет класс объектов, состоящий из всех  гипергрупп над фиксированной группой $H$,
совпадающей с мультипликативной группой поля $k$,
и с системой структурных отображений $\Omega = (\Phi, \Psi, \Xi, \Lambda)$ такой что

(i) структурные отображения  $\Psi$ и $\Lambda$ тривиальны,

(ii) $(M, \Xi)$ является абелевой группой,

(iii)  структурное отображение $\Phi$ и  группа $H$ подчиняются следующим двум условиям.

\hspace{10mm}
Во-первых, группа $H$ абелева.
Второе условие формулируется сложнее.
Сначала заметим, что  в силу условия (ii)
на множестве $End (M, \Xi)$ эндоморфизмов (абелевой) группы $(M, \Xi)$
можно определить  операцию сложения по формуле
$$
(\varphi + \psi) (a) = \varphi (a) + \psi (a) := \Xi (\varphi (a). \psi (a)), \quad \varphi, \psi \in End (M, \Xi), \, a \in M,
$$
причем множество $End (M, \Xi)$ вместе с этой операцией становится абелевой группой.
Нейтральным элементом $\theta$ этой группы будет "`нулевой эндоморфизм"',
переводящий любой  элемент $a \in M$ в нейтральный элемент $o$ группы $(M, \Xi)$.
В силу соотношения $(A1)$ определения гипергруппы над группой
существует каноническое инъективное отображение
$t: H \rightarrow End (M, \Xi)$.
Второе условие заключается в том, что множество
$$
(F5)  \hspace{45mm}  k = t(H) \cup \{ \theta \}  \hspace{50mm}
$$
должно быть замкнутым относительно вышеопределенной операции сложения в $End (M, \Xi)$, и
эта операция сложения вместе с операцией умножения, определяемой по правилу
$$
t(\alpha) \cdot t(\beta) = t(\alpha \cdot \beta) \, \, (\alpha, \beta \in H),  \quad
\theta \cdot \varphi = \theta =  \varphi \cdot \theta \, \, (\varphi \in k),
$$
должны задавать на $k$ структуру поля.

\hspace{10mm}
Обратный к функтору $S(L)$ - это функтор из $im S(L)$ в ${\mathcal L}_k$,
который сопоставляет гипергруппе над группой $M_H$ из категории
$im S(L)$ с системой структурных отображений $\Omega = (\Phi, \Psi, \Xi, \Lambda)$
линейное пространство $L$ над полем $k$, где
поле $k$ задается как множество $(F5)$ вместе с определенными на этом множестве
операциями сложения и умножения,
а линейное пространство  $L$ над полем $k$ есть
множество $M$ вместе с  операцией сложения $\Xi$
и операцией умножения на элементы из поля $k$, совпадающей
с действием соответствующих элементов из $End (M, \Xi) \supset k$.


\vspace{5mm}

\hspace{10mm}
{\bf 10.}
{\bf Предложение.}
{\it Категория гипергрупп над группой $\mathcal{H}g$ содержит подкатегорию,
изоморфную категории полей $\mathcal{F}$.}

\hspace{10mm}
{\it Доказательство.}
Конструкция функтора $S(F): {\mathcal F} \rightarrow  {\mathcal H}g$ аналогична конструкции функтора $S(L)$.
Каждому полю $F$ сопоставляется гипергруппа над группой $M_H$ с системой структурных отображений
$\Omega = (\Phi, \Psi, \Xi, \Lambda)$, где

$M$ совпадает с базовым множеством поля $F$,

$H$ - мультипликативная группа поля $F$,

$\Phi$ является ограничением операции умножения поля $F$,

$\Xi$ - операция сложения поля $F$,

$\Psi$ и $\Lambda$ - тривиальные структурные отображения.

На подкатегорию $im S(F)$ накладывается дополнительное
к накладываемым на подкатегорию $im S(L)$ условие,
заключающееся в том, что
аддитивная группа поля $k$ должна быть канонически изоморфной абелевой группе $(M, \Xi)$.


\vspace{5mm}

\hspace{10mm}
{\bf 11.}
{\bf Замечание.}
Если отказаться от условия  коммутативности операции умножения в поле и,
соответственно, от условия абелевости группы $H$ в (iii) (см. доказательство предложения 9),
то получим доказательство утверждения, что
категория гипергрупп над группой содержит в качестве подкатегорий
категорию линейных пространств над телом и категорию тел.
Подчеркнем, однако, что согласно известному результату Веддерберна
это замечание ничего нового не прибавляет в случае конечных тел.

\vspace{5mm}

\hspace{10mm}
{\bf 12.}
{\bf Замечание.}
Гипергруппы над группой уже успешно применены в [6] дкя обобщения
теоремы Шрайера о расширениях группы с помощью коммутативной группы
([1], стр. 185, Теорема 7.34) и для описания полупрямых произведений групп.

\end{document}